\documentclass{tran-l}
\usepackage{amssymb,latexsym, amscd}
\usepackage[all]{xy}
\usepackage{graphicx}
 \newtheorem{theorem}{Theorem}[subsection]
 \newtheorem{cor}[theorem]{Corollary}
 \newtheorem{lemma}[theorem]{Lemma}
 \newtheorem{proposition}[theorem]{Proposition}
 \theoremstyle{definition}
 \newtheorem{definition}[theorem]{Definition}
 \theoremstyle{definition}
 \newtheorem{example}[theorem]{Example}
 \theoremstyle{remark}
 \newtheorem{rem}[theorem]{Remark}
 \numberwithin{equation}{subsection}

\usepackage{latexsym}
\usepackage{amssymb}

\newcommand{\ben}{\begin{equation}}
\newcommand{\een}{\end{equation}}

\newcommand{\integer}{\ensuremath{{\mathbb Z}}}
\newcommand{\naturals}{\ensuremath{{\mathbb N}}}
\newcommand{\real}{\ensuremath{{\mathbb R}}}
\newcommand{\complex}{\ensuremath{{\mathbb C}}}

\newcommand{\rational}{\ensuremath{{\mathbb Q}}}

\newcommand{\XX}{{\mathcal X}}

\newcommand{\UU}{{\mathcal U}}

\newcommand{\GG}{{\mathcal G}}
\newcommand{\CC}{\mathcal{C}}

\newcommand{\LL}{\mathcal{L}}

\newcommand{\OO}{\mathcal{O}}
\newcommand{\YY}{\mathcal{Y}}
\newcommand{\TT}{\mathcal{T}}

\newcommand{\Hom}{\mathrm{Hom}}

\newcommand{\Xx}{\mathsf{X}}
\newcommand{\Gg}{\mathsf{G}}
\newcommand{\Hh}{\mathsf{H}}

\newcommand{\Ss}{\mathsf{S}}

\newcommand{\target}{\mathsf{t}}
\newcommand{\source}{\mathsf{s}}
\newcommand{\ident}{\mathsf{e}}
\newcommand{\invers}{\mathsf{i}}
\newcommand{\mult}{\mathsf{m}}

\newcommand{\To}{\longrightarrow}

\newcommand{\timests}{\: {}_{\target}  \! \times_{\source}}

\newcommand{\Vector}[1]{\stackrel{\rightarrow}{#1}}
\newcommand{\smallloop}{{\LL_s}}
\newcommand{\coarse}[1] {{#1} / \! \! \sim}
\newcommand{\Global}[2]{\, {#1} \! \! \rtimes  \! {#2}}

\newcommand{\Conj}{\ensuremath{{\mathrm{Conj}}}}

\begin{document}

\title{Inertia Orbifolds, Configuration Spaces and the Ghost Loop Space.}
\author{ Ernesto Lupercio and Bernardo Uribe}
\thanks{The first author was partially supported by the National Science
Foundation}

\address{Department of Mathematics, University of Wisconsin at Madison, Madison, WI 53706}
\address{Max-Planck-Institut f\"{u}r Mathematik, Vivatsgasse
7, D-53111 Bonn, Germany Postal Address: PO.Box: 7280, D-53072
Bonn} \email{ lupercio@math.wisc.edu \\ uribe@mpim-bonn.mpg.de}

\begin{abstract}
In this paper we define and study the \emph{ghost loop orbifold}
${\mathcal{L}_s}{ \mathsf{X}}$ of an orbifold $\mathsf{X}$
consisting of those loops that remain constant in the coarse
moduli space of $\mathsf{X}$. We construct a configuration space
model for ${\mathcal{L}_s}{ \mathsf{X}}$ using an idea of G.
Segal. From this we exhibit the relation between the Hochschild and
cyclic homologies of the inertia orbifold of ${ \mathsf{X}}$
(that generate the so-called twisted sectors in string theory)
and the ordinary and equivariant homologies of ${\mathcal{L}_s}{
\mathsf{X}}$. We also show how this clarifies the relation
between orbifold K-theory, Chen-Ruan orbifold cohomology,
Hochschild homology, and periodic cyclic homology.
\end{abstract}

\maketitle

\section{Introduction}

Given an orbifold $X$, Kawasaki \cite{Kawasaki} has defined an associated orbifold that is often called the 
\emph{twisted sector} orbifold $\wedge \Xx$  \cite{Ruan}. This orbifold is also known as the 
\emph{inertia orbifold} associated to $\Xx$. Essentially, this orbifold is made of pairs $(x,(g))$ where
 $x \in X$ and $(g) \in \Conj(G_x)$ where $G_x$ is the stabilizer of $\Xx$ at $x$. A full description of
  $\wedge \Xx$ needs, of course, a precise definition of the orbifold structure on the space of such pairs.
   This new orbifold appears naturally from the point of view of many geometric questions, and again,
   from the perspective of string theory its introduction is well motivated.  For example, the orbifold
   $K$-theory of $\Xx$ is rationally isomorphic to the cohomology of $\wedge \Xx$ by the Chern character.
   The purpose of this paper is to study the topology of this orbifold using standard algebraic topology.
    In particular we will not use the theory of $C^*$-algebras.

We should clarify here what we mean by the topology of an orbifold. An orbifold $\Xx$, let us recall,
 is a pair $(X,\UU)$ where $X$ is a topological space and $\UU$ is an equivalence class of orbifold atlases.
  An orbifold atlas $\UU=(\pi_i, U_i, G_i, V_i)_i$ is made up of orbifold charts so that
   $V_i$ is an open subset of $X$, $U_i$ is a manifold, $G_i$ acts on $U_i$ via diffeomorphisms and
   $\pi_i\colon U_i \to V_i$ is a quotient map that induces a homeomorphism $U_i/G_i \simeq V_i$;
      it is important to recall that all the groups $G_i$ are finite. The gluing conditions are subtle.
       At a first level of approximation  the topology of an orbifold $\Xx$ is simply the topology of
       the topological space $X$. This naive approach is neglecting altogether the orbifold structure
        $\UU$. To  consider the orbifold structure from the point of view of algebraic topology we
	follow Haefliger \cite{Haefliger} and  Moerdijk-Pronk \cite{MoerdijkPronk}, and associate a
	topological category $\Gg$ to the orbifold $\Xx$ (actually the category is a groupoid.) Then
	we consider the classifying space $B\Gg$. As Moerdijk \cite{MoerdijkTopos} has shown the homotopy
	 type of the orbifold $\Xx$ (seen as a topological stack or as a topos, in the sense of \'{e}tale
	  homotopy theory) is exactly the same thing as the homotopy type of $B\Gg$. The situation is
	  analogous to that of studying the space $BG$ associated to a finite group $G$. In fact this
	  situation inspires much of what follows.

The category $\Gg$ associated to an orbifold is \'{e}tale (the structure maps are local diffeomorphisms)
and \emph{stable}  (the sets $\Hom_\Gg (x,y)$ are finite.) But in fact the constructions of this paper
 are valid for a general topological category $\Gg$. Given a groupoid $\Gg$ its inertia groupoid
 $\wedge \Gg$ is again a groupoid. Its set of objects is $\Hom(\integer, \Gg)$ (here we mean all
 functors from the category with one object and a $\integer$ worth of arrows, to the category $\Gg$.)
  The morphisms of $\wedge \Gg$ are the obvious ones (induced by the action of $\Gg$ on
   $\Hom(\integer, \Gg)$.) As we said before the orbifold $K$-theory of  $\Xx$ computes the cohomology
   of $\wedge \Xx$ rationally.

For every functor $\integer \to \Gg$ there is a corresponding map $S^1=B\integer \to B\Gg$, namely
 an element in $\LL B \Gg:= {\rm Map}(S^1,B\Gg)$. Not every element of $\LL B \Gg$ comes from such
  a functor. We define then the space of \emph{ghost loops} $\smallloop B\Gg$ of $B\Gg$ as the subspace
   of $\LL B \Gg$ consisting of loops so that the composition $S^1 \to B \Gg \to X$ is constant
    (Here $B \Gg \to X$ is the natural quotient projection.) In section \ref{sectionghostloopspace} below
    we prove the following theorem.

\begin{theorem} If $\Gg$ is an orbifold groupoid then the
 natural map $\CC_{S^1}(\Gg) \longrightarrow \smallloop B \Gg $ induces a weak homotopy equivalence as well as a 
weak $S^1$-homotopy equivalence.
\end{theorem}

The main idea in the proof of this theorem consist in the introduction of a \emph{configuration space model}
 for the space $B\wedge \Gg$. Following a beautiful idea of Segal and Burghelea in the case of a group, we
  define the space $\CC_{S^1}(\Gg)$ to be the configuration space of particles in $S^1$ labeled by morphisms
   of $\Gg$ (so that every morphism composes with the ones adjacent to it.)

\begin{theorem} For $\Gg$ a topological groupoid there
 is a natural homeomorphism $B \wedge \Gg \cong \CC_{S^1}( \Gg)$.
 \end{theorem}

For a finite group there is a well known relation between the homology of $\LL B G$ and the Hochschild and
cyclic homologies of $R[G]$ \cite{ BurgheleaFiedorowicz, Goodwillie, Jones}. Using the previous theorem we
 prove the following result, that can be seen as a generalization of this relation to the case of
 smooth orbifolds.

\begin{theorem}
Let $\Gg$ be an orbifold groupoid and $\wedge \Gg$ its inertia groupoid. Then there are
canonical isomorphisms
\begin{eqnarray*}
HH_*((\wedge \Gg_*,t_*)) \stackrel{\cong}{\To} H_*(B \wedge \Gg) & \stackrel{\cong}{\To} & H_*(\CC_{S^1}(\Gg))
\stackrel{\cong}{\To} H_*(\smallloop B \Gg)\\
  HC_*((\wedge \Gg_*,t_*)) \stackrel{\cong}{\To} H^{S^1}_*(B \wedge \Gg) & \stackrel{\cong}{\To} & H^{S^1}_*(\CC_{S^1}(\Gg))
\stackrel{\cong}{\To} H^{S^1}_*(\smallloop B \Gg)
\end{eqnarray*}
\end{theorem}

We should point out that these results are related to the theorems of Brylinski-Nistor \cite{BrylinskiNistor}
 and Crainic \cite{Crainic} obtained using the theory of $C^*$-algebras.

Finally as a direct consequence of the previous theorem we show the following isomorphisms involving the
 Chen-Ruan cohomology of the orbifold
$H_{\mathrm{orb}}^*(\Gg ; \complex)$ and also its orbifold $K$-theory \cite{ AdemRuan, ChenRuan, LupercioUribe1}.

\begin{theorem}
For $\Gg$ a compact complex $SL$-orbifold the following holds,
$$\prod_{m \in n + 2 \integer} HH_m((\wedge \Gg_*, t_*) ; \complex)
\cong \prod_{m \in n + 2 \integer}H_{\mathrm{orb}}^m(\Gg ; \complex)$$
and if $\Gg$ is reduced
$$K_{\mathrm{orb}}^n(\Xx)\otimes \complex
\cong \prod_{m \in n + 2 \integer} H_{\mathrm{orb}}^m(\Xx, \complex)\cong HP_n(\wedge \Gg, \complex).$$
\end{theorem}

We would like to call the attention of the reader to the recent papers of
Baranovsky \cite{Baranovsky}, and C\u ald\u araru-Giaquinto-Witherspoon \cite{CaldararuGiaquintoWitherspoon}
 where relations
between orbifold cohomology and Hochschild and periodic cyclic homology are
also explored.

We would like to thank conversations with A. Adem, D. Berenstein, M. Karoubi,  B. Oliver, J. Robbin,
and Y. Ruan and especially we would like to thank I. Moerdijk for helping us with the proof of
theorem \ref{w.h.e.conf=ghost}.

We dedicate this paper to Prof. Graeme B. Segal on the occasion of his $60^{\mathrm{th}}$
birthday.
\section{Groupoids}
The groupoids we will consider are small categories $\Gg$ in which every morphism
 is invertible. By $\Gg_1$ and $\Gg_0$ we will denote the space of morphisms (arrows)
and of objects respectively, and the structure maps

        $$\xymatrix{
        \Gg_1 \timests \Gg_1 \ar[r]^{\mult} & \Gg_1 \ar[r]^{\invers} &
        \Gg_1 \ar@<.5ex>[r]^{\source} \ar@<-.5ex>[r]_{\target} & \Gg_0 \ar[r]^{\ident} & \Gg_1
        }$$
where $\source$ and $\target$ are the source and the target maps, $\mult$ is the composition of two arrows,
$\invers$ is the inverse and $\ident$ gives the identity arrow over every object.

 The groupoid will be called {\it topological (smooth)} if the sets $\Gg_1$ and $\Gg_0$
 and the structure maps belong to the category of topological spaces (smooth manifolds).
In the case of a smooth groupoid we
will also require that the maps $\source$ and $\target$ must be submersions, so that
$\Gg_1 \timests \Gg_1 $ is also a manifold (smooth groupoids are also known by the name
of {\it Lie groupoids} \cite{Moerdijk2002}).

A topological (smooth) groupoid is called {\it \'{e}tale} if all the structure maps
are local homeomorphisms (local diffeomorphisms). For an \'{e}tale groupoid we will
mean a topological \'{e}tale groupoid. In what follows, sometimes the kind of groupoid
will not be specified, but it will be clear from the context to which one we are referring to.
We will always denote groupoids by letters of the type $\Gg,\Hh,\Ss$. Orbifolds are an special kind of 
\'{e}tale groupoids, they have the peculiarity
the inverse imagine of a point under the map $(\source,\target): \Gg_1 \to \Gg_0\times \Gg_0$ is always
finite. Groupoids with this property are called \emph{stable}\footnote{We owe this terminology to J. Robbin.
 It is inspired in the stable map compactification for holomorphic maps.}. Whenever we write orbifold, 
 a stable \'{e}tale smooth groupoid will be understood.

A morphism of groupoids $\Psi: \Hh \to \Gg$ is a pair of maps
$\Psi_i: \Hh_i \to \Gg_i$ $i=1,2$ such that they commute with the
structure maps. The maps $\Psi_i$ will be continuous (smooth)
depending on which category we are working on.

For a groupoid $\Gg$, we denote by $\Gg_n$ the space of $n$-arrows
$$x_0 \stackrel{g_1}{\to} x_1 \stackrel{g_2}{\to} \cdots \stackrel{g_n}{\to} x_n.$$
The spaces $\Gg_n$ $(n \geq 0)$ form a simplicial space $(\Gg_*, d^i_n, s^i_n)$
\begin{eqnarray}        
\xymatrix{
        \cdots \ar@<1.5ex>[r] \ar@<.5ex>[r] \ar@<-.5ex>[r] \ar@<-1.5ex>[r] &
        \Gg_2 \ar@<1ex>[r] \ar[r] \ar@<-1ex>[r] &
        \Gg_1 \ar@<.5ex>[r]^\source \ar@<-.5ex>[r]_\target & \Gg_0
        } 
\end{eqnarray}
together with the face $d_n^i: \Gg_n \to \Gg_{n-1}$ and degeneracy maps $s_n^i : \Gg_n \to \Gg_{n+1}$ (see \ref{nervegroupoid})
form what is called the {\it nerve of the groupoid}. Its geometric realization $|\Gg_*|$ (see Appendix)
 is the classifying space of $\Gg$, also denoted by $B\Gg$.

\section{Configuration space  models for groupoids.}

In this section we will apply ideas of G. Segal  
regarding the configuration
spaces on groups to the more general setting of topological groupoids. On what follows
the groupoids will be topological unless otherwise stated. 

By $S^1 = \real/ \integer$ we will denote the circle of radius $\frac{1}{2 \pi}$ and by $R$ we will
denote the open interval $(0,1)$. Let $\iota : R \to S^1$ be the composition of the projection and
inclusion maps $(0,1) \hookrightarrow \real \to \real/ \integer$. Let $A$ be either $R$ or $S^1$.

For $\alpha : A \to \Gg_1$ a map, we will call the support of $\alpha$ the set:
$${\rm supp} (\alpha) := \{ a \in A | \alpha (a) \notin \Gg_0 \}$$ 
where we consider the set $\Gg_0$ as the image under the identity map
of $\ident(\Gg_0) \subset \Gg_1$.

\begin{definition} \label{definitionconfiguration}
A configuration in $A$ with values in $\Gg$ is a map $\alpha : A \to \Gg_1$
with finite support
$${\rm supp} (\alpha) := \{x_1, x_2, \dots, x_n \} \ \ \ \ \  0 \leq x_i < x_{i+1} <1$$
such that for $g_i := \alpha(x_i)$
$$\target (g_i) = \source (g_{i+1}) \ \ \ \ \ \ \rm{for}\ \  1\leq i <n$$
and when $A=S^1$
$$\source (g_1) = \target (g_n).$$

\begin{center}
\bigskip

\includegraphics[height=2.2in]{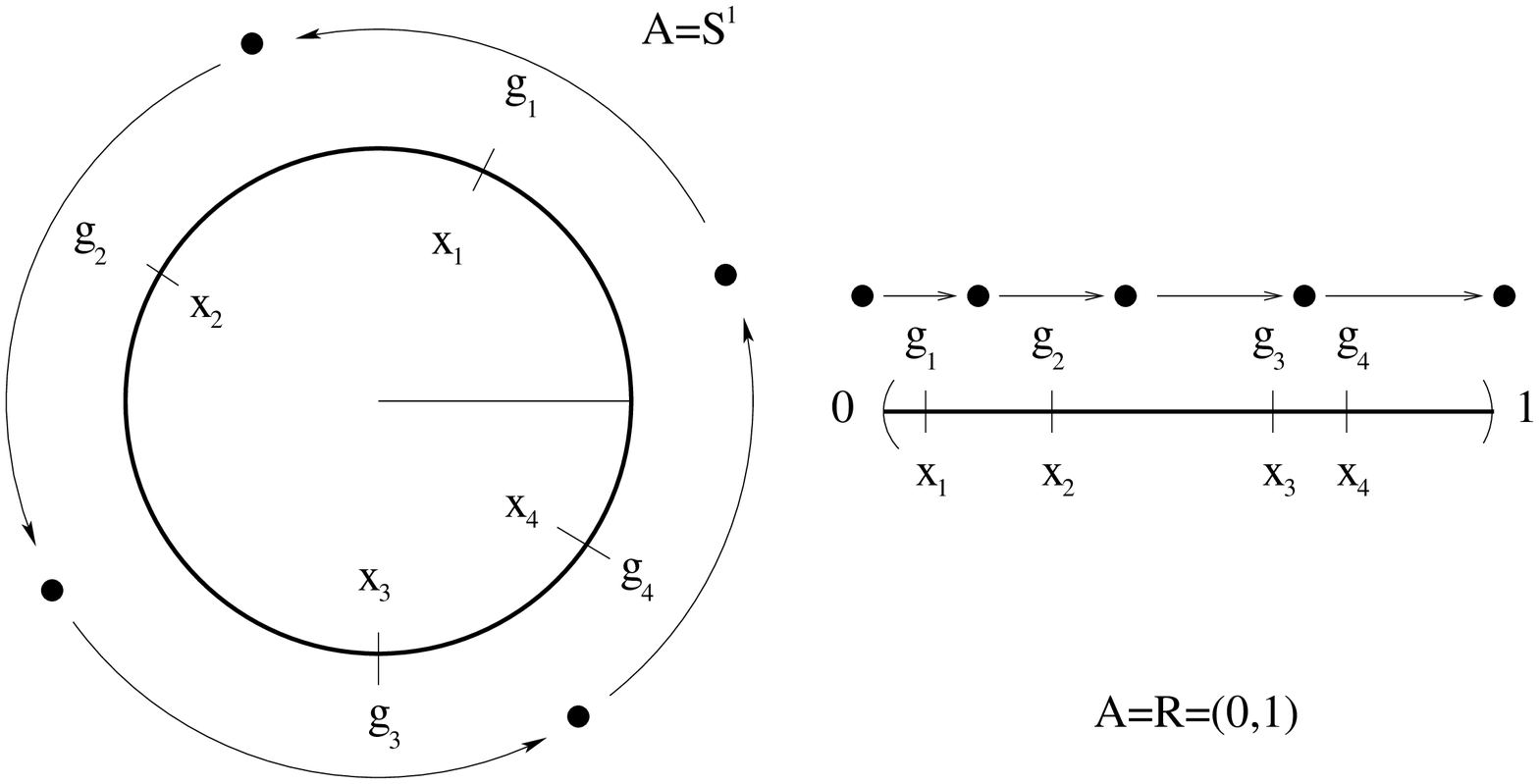}
\end{center}
Also we require that the image of $\alpha$ outside its support be determined by the $g_i$'s
 as follows, for $x_i < y < x_{i+1}$
$$\alpha(y) = \ident \target (g_i),$$
$\alpha(y)= \ident \source (g_{1})$ for $y < x_1$ and $\alpha(y)= \ident \target (g_{n})$ for $y > x_n$.

The collection of configurations will be denoted by $\CC_A(\Gg)$.
\end{definition}

We are going to endow these configuration spaces with natural topologies in such a way that
when two morphisms of the groupoid collapse they get composed, and when a morphism gets sent
to infinity it disappears (in the case $A=R$).

A base for the topology of $\CC_A(\Gg)$ is given by the sets $\OO(U_1,U_2, \dots, U_k ; I_1, I_2, \dots , I_k)$
where $U_i$ are open sets $\Gg_1$ and $I_i$ are connected, open and disjoint intervals (on the circle
the orientation of $S^1$ provides us with a total order in the $I_i$'s). The set
$\OO(U_1, \dots, U_k ; I_1, \dots , I_k)$ denotes the set of configurations $\alpha$ whose support
consists of the points $x_i^j$ with
$$x_1^j < x_2^j < \cdots < x_{l_j}^j,  \ \ \ \ \ \ \ x_i^j \in I_j \ \ \ \ \
 {\rm and} \ \ \ \ \ g_1^j g_2^j \cdots g_{l_j}^j \in U_j$$
for $g_i^j := \alpha(x_i^j)$.

There are natural maps relating these configuration spaces:
\begin{itemize}
\item by rotation of a configuration by an  angle  $\theta$, 
\begin{eqnarray}
\rho : S^1 \times \CC_{S^1}(\Gg)  \to \CC_{S^1}(\Gg) \label{rotationmap}
\end{eqnarray}
where $\rho(\theta,\alpha)(x) := \alpha(\theta + x)$,
\begin{center}
\includegraphics[height=1.5in]{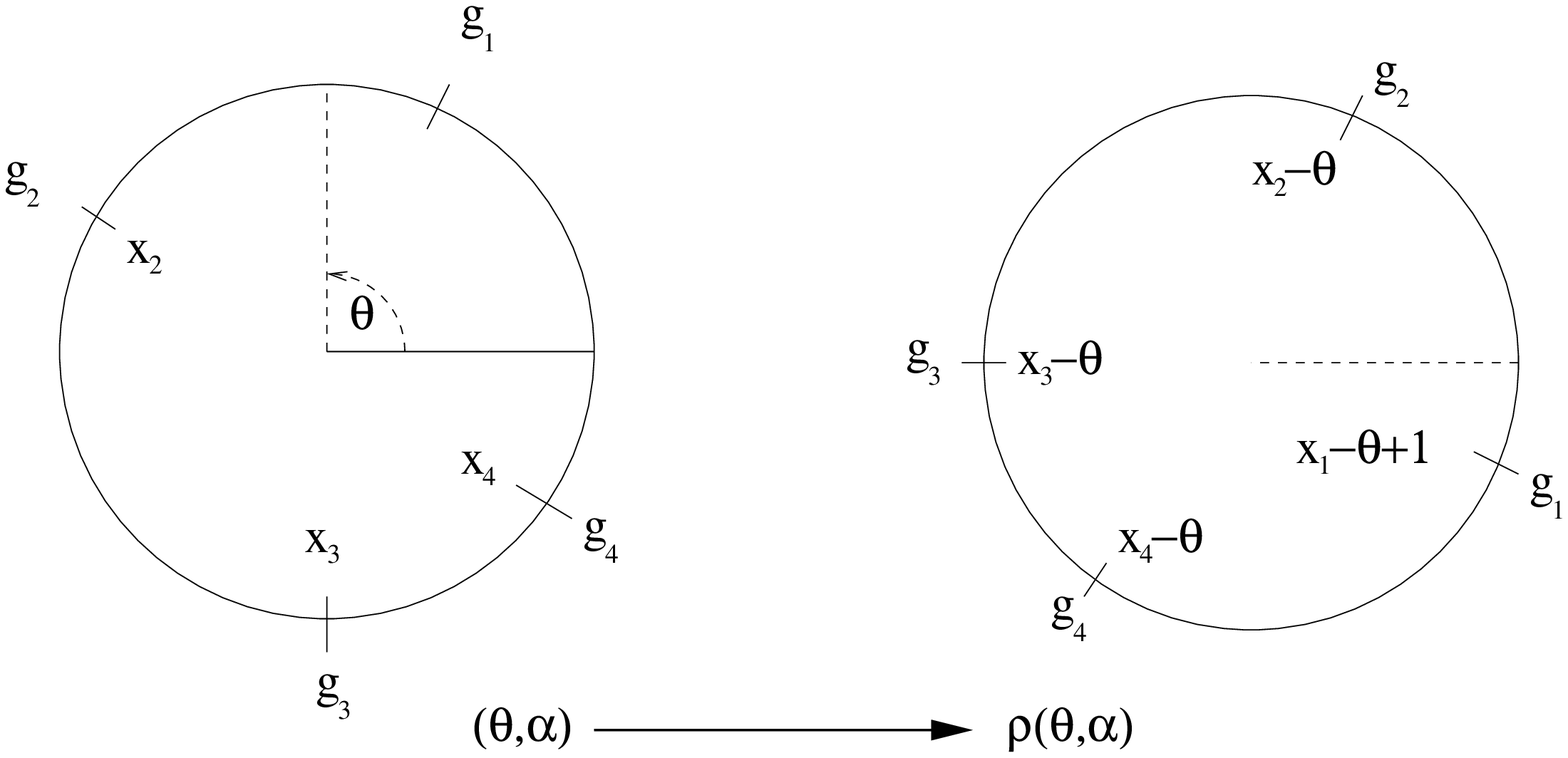}
\end{center}
\item by sending to infinity the morphism at the angle $\theta$,
\begin{eqnarray}
\varepsilon: \CC_{S^1}(\Gg) \times S^1\to \CC_R(\Gg) \label{evaluationmap}
\end{eqnarray}
where $\varepsilon(\alpha, \theta) \in \CC_R(\Gg)$ such that for $x \in (0,1)$,
 $\varepsilon(\alpha,\theta)(x) := \alpha(\theta + \iota(x))$
\begin{center}
\includegraphics[height=1.5in]{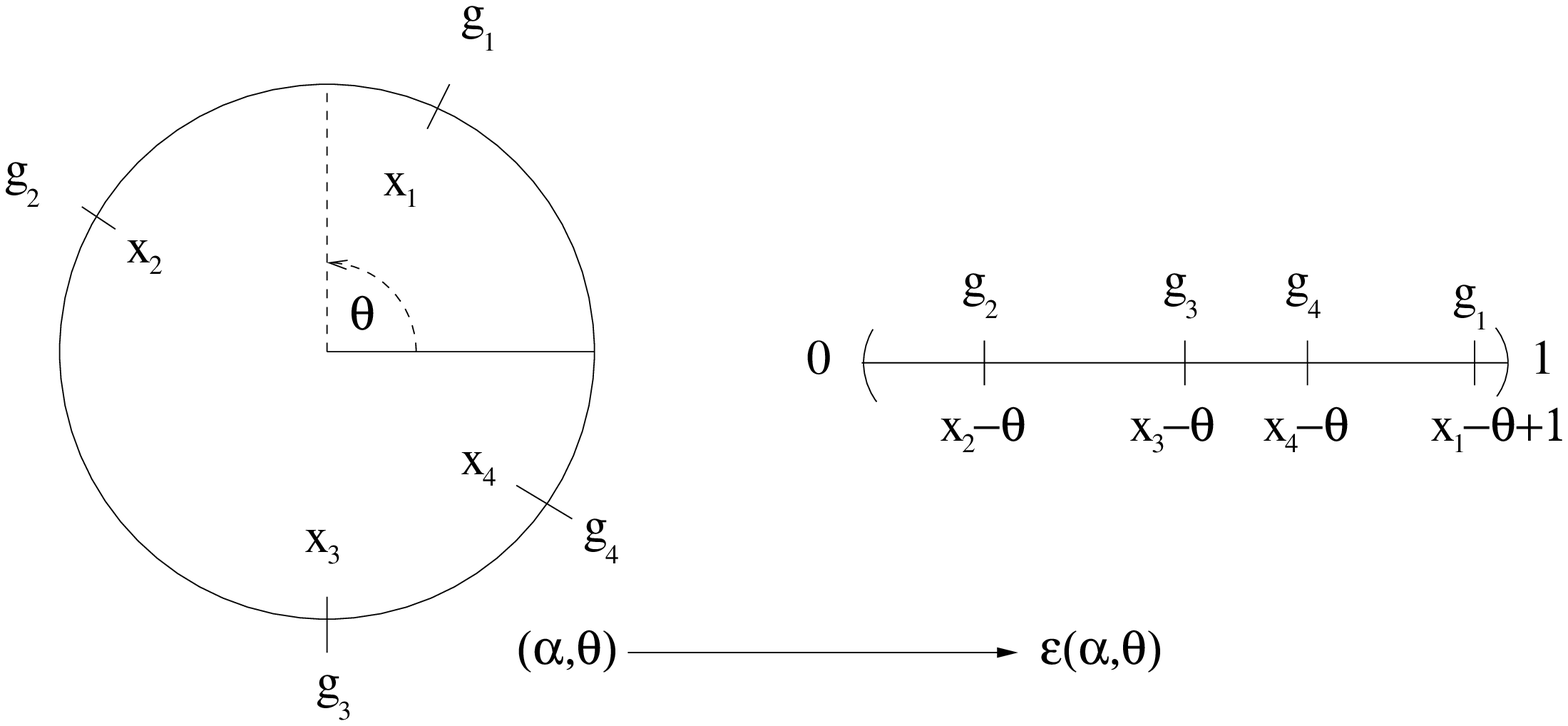}

\includegraphics[height=1.5in]{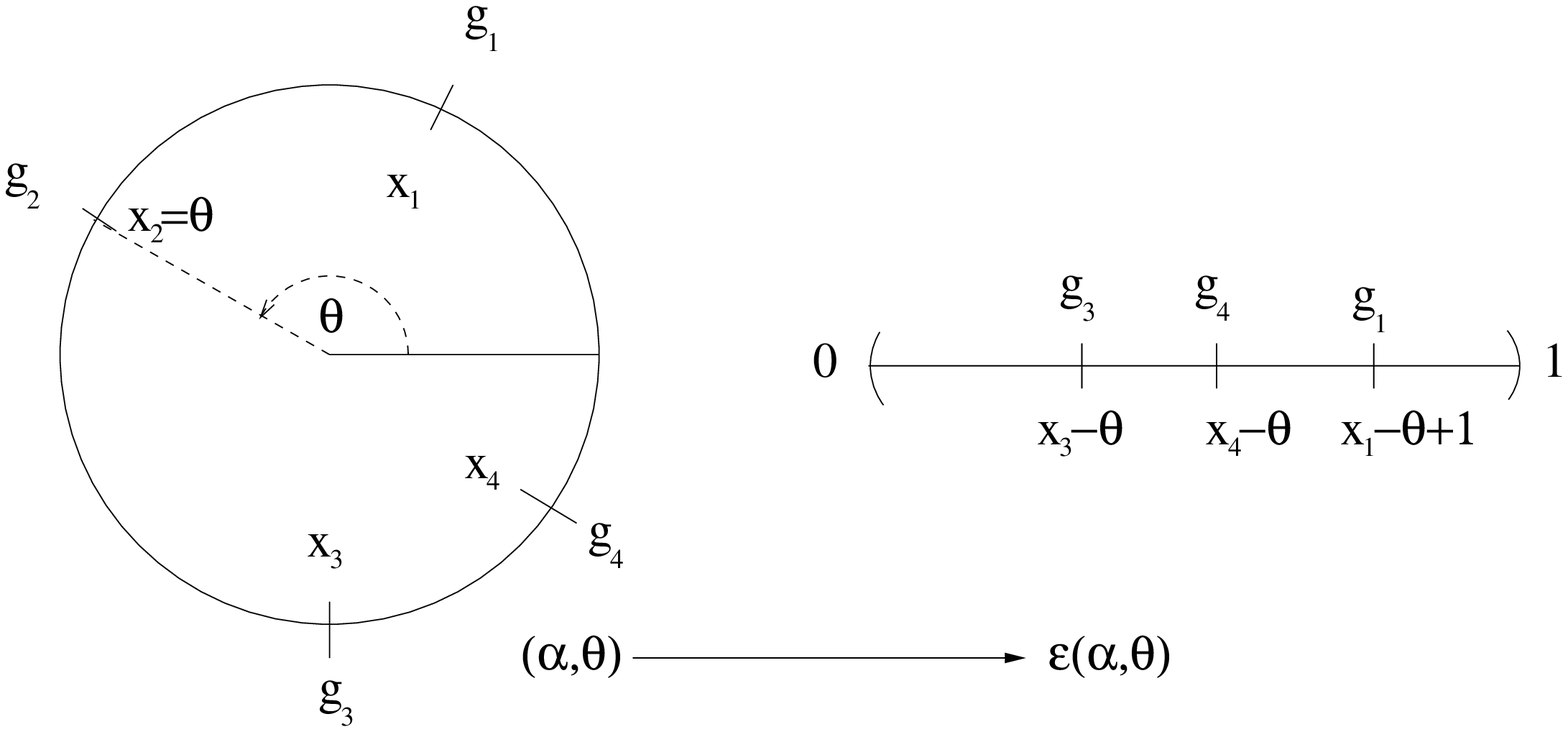}
\end{center}
\end{itemize} 

With the topology previously defined it is clear that the maps $\rho$ and $\varepsilon$ are continuous.
The set $\CC_A(\Gg)$ can be filtered
$$* = \CC_A^0(\Gg) \subset \CC_A^1(\Gg) \cdots \CC_A^n(\Gg) \subset \CC_A^{n+1}(\Gg) \subset \cdots$$
where $\CC_A^n(\Gg):=\{\alpha \in \CC_A(\Gg) \ : \ |{\rm supp} (\alpha)| \leq n \}$. It is also
clear that the maps $\rho$ and $\varepsilon$ are compatible with the filtration.

\begin{proposition} \label{propositionBG}
$\CC_R(\Gg)$ is canonically isomorphic to $B\Gg$.
\end{proposition}
\begin{proof}
First of all let's recall how the classifying space of a groupoid is built. For a 
groupoid $\Gg$ we associate a simplicial space $\Gg_*=(\Gg_n, d_n^i, s_n^i)$
where
\begin{eqnarray}
\Gg_n := \underbrace{\Gg_1 \timests \Gg_1 \timests \cdots \timests \Gg_1}_{n} \label{nervegroupoid}
\end{eqnarray}
$$d_n^i(g_1, g_2, \dots, g_n) = \left\{
\begin{array}{cc}
(g_2, \dots, g_n) & \ i=0 \\
(g_1, \dots, g_{i-1}, g_ig_{i+1}, g_{i+2}, \dots, g_n) & i \neq 0,n \\
(g_1, \dots, g_{n-1}) & i=n
\end{array} \right. $$
$$s_n^i(g_1, \dots, g_n) = (g_1, \dots, g_{i-1}, \ident\source(g_{i}), g_i , \dots, g_n).$$
The geometric realization $| \Gg_*|$ of the simplicial space $\Gg_*$ is what
we denote $B\Gg$, the classifying space of the groupoid. This geometric realization
is defined as
$$B \Gg=|\Gg_*|:= \bigsqcup_{n=0}^\infty \Gg_n \times \Delta_n \ / \sim = 
\bigsqcup_{n=0}^\infty (\Gg_n- {\rm Deg} \ \Gg_n) \times \Delta_n \ / \sim$$
with ${\rm Deg} \ \Gg_n = \bigcup_{i=0}^{n-1} s_{n-1}^i(\Gg_{n-1})$ and
 $\sim$ the equivalence relation generated by the identifications
$$(s^i_{n-1}(x), \Vector{t}) \sim (x, \sigma_n^i(\Vector{t}) \ \ \ {\rm and}
\ \ \ (d_{n+1}^i(y), \Vector{t}) \sim (y, \delta_n^i(\Vector{t})$$
with $x \in \Gg_{n-1}$, $y \in \Gg_{n+1}$, $\Vector{t} \in \Delta_n$ and
$\Delta^*=(\Delta_n, \delta_n^i, \sigma_n^i)$ the cosimplicial set of the $n$-simplices
with their face and degeneracy maps:
$$\Delta_n := \{\Vector{t}=(t_0, \dots, t_n) \in \real^{n+1} \}$$
$$\delta_n^i(t_0, \dots, t_n) := (t_0, \dots, t_i, 0, t_{i+1} \dots, t_n)$$
$$\sigma_n^i(t_0, \dots, t_n) := (t_0, \dots, t_{i-1}, t_i +t_{i+1}, t_{i+2}, \dots, t_n).$$

Taking out the images of the face maps, in other words, only considering the
set
$$ {\rm Int} \ \Delta_n := \Delta_n - \bigcup_{i=0}^{n-1}\delta_{n-1}^i(\Delta_{n-1})$$
then it is clear that as a set 
 $$B\Gg = \bigsqcup (\Gg_n - {\rm Deg} \ \Gg_n) \times {\rm Int} \ \Delta_n.$$

Now we can construct a bijective correspondence between the points in $\CC_R(\Gg)$ and the
points in $B\Gg$ in the following way: to the configuration $\alpha : R \to \Gg_1$
with ${\rm supp}(\alpha) = \{x_1, x_2, \dots, x_n \}$ and $\alpha(x_i) = g_i$ we associate
the point 
$$(g_1, g_2, \dots, g_n; t_0, t_1, \dots, t_n) \in 
(\Gg_n - {\rm Deg} \ \Gg_n) \times {\rm Int} \ \Delta_n$$
with $t_0 = x_1$, $t_1= x_2 - x_1$, \dots $t_n = 1 - x_n$. It is easy to see now that this
identification is a homeomorphism. 
\end{proof}

\section{The Inertia Groupoid}\label{Inertiagroupoid}

The Inertia groupoid $\wedge \Gg$ of a groupoid $\Gg$ is defined as follows.

 An object $a$  ($a \in \wedge \Gg_0$) is an arrow of $\Gg_1$
such that its source and target are equal, so
$$\wedge\Gg_0 = \{ a \in \Gg_1 | \source(a)=\target(a)\}$$
and a morphism $v \in \wedge \Gg_1$ joining two objects $a,b$
         $$\xymatrix{
         \circ \ar@(ul,ur)[]^a \ar@/^/[r]^v & \circ \ar@(ul,ur)[]^b}$$
is an arrow in $\Gg_1$ such that $a\cdot v = v\cdot b$; $b$ can be seen as
the conjugate of $a$ by $v$ because $b = v^{-1} \cdot a \cdot v$. But in order
to identify the arrow $v$ we need to keep track of its source $a$, then we can consider
the morphisms of $\wedge \Gg$ as
$$\wedge \Gg_1 = \{ (a,v) \in \Gg_2 | a \in \wedge \Gg_0 \}$$
where $\source(a,v) = a$ and $\target(a,v) = v^{-1} \cdot a\cdot v$.

Let $\wedge \Gg_* = (\wedge \Gg_n , d_n^i, s_n^i)$ be the simplicial space associated
to the inertia groupoid, where
$$\wedge \Gg_n =  \{ (a,v_1,v_2, \dots, v_n) \in \Gg_{n+1} | a \in \wedge \Gg_0 \}$$
$$d_n^i(a,v_1,v_2, \dots, v_n) = \left\{
\begin{array}{cc}
(v_1^{-1}av_1, v_2,\dots,v_n) & i = 0 \\
(a,v_1, \dots, v_{i-1}v_{i}, \dots, v_n) & i \neq 0,n \\
(a,v_1, \dots, v_{n-1}) & i=n
\end{array} \right.$$
$$s^i_n(a,v_1,v_2, \dots, v_n) = (a, v_1, \dots, v_{i-1}, \ident\source (v_i), v_i, \dots, v_n)$$

It turns out that there is another simplicial space $\YY_*(\Gg)$ associated to a groupoid
which is isomorphic to $\wedge \Gg_*$ and that is related in a more direct way to the
configuration space $\CC_{S^1}(\Gg)$.
$\YY_*(\Gg)$ is defined as follows:
$$\YY_n(\Gg) := \{ (g_0, g_1, \dots, g_n) \in \Gg_n \ |\ \source(g_0) = \target(g_n) \}$$
$$d_n^i(g_0, g_1, \dots, g_n) = \left\{
\begin{array}{cc}
(g_0, \dots, g_{i-1},g_ig_{i+1}, g_{i+2} , \dots, g_n) & i \leq n-1 \\
(g_ng_0,g_1, \dots , g_{n-1}) & i=n 
\end{array} \right.$$
$$s_n^i(g_0, g_1, \dots, g_n) = (g_0, \dots, g_i, \ident \source(g_{i+1}), g_{i+1}, \dots,
g_n)$$

\begin{proposition}
The simplicial spaces $\wedge \Gg_*$ and $\YY_*(\Gg)$ are isomorphic.
\end{proposition}
\begin{proof}
Define the maps
\begin{eqnarray*}
f_n : \wedge \Gg_n & \to & \YY_n(\Gg) \\
(a,v_1,v_2, \dots, v_n) & \mapsto & (v_n^{-1}\cdots v_2^{-1}v_1^{-1}a,v_1,v_2, \dots, v_n)
\end{eqnarray*}

\begin{eqnarray*}
h_n : \YY_n(\Gg) & \to & \wedge \Gg_n\\
(g_0,g_1, \dots, g_n) & \mapsto & (g_1g_2\cdots g_ng_0,g_1,g_2, \dots, g_n),
\end{eqnarray*}
clearly $h_n \circ f_n = {\rm id}_{\wedge \Gg_n}$ and 
$f_n \circ h_n = {\rm id}_{\YY_n(\Gg)}$. It is easy to see that
the maps $f_n$ and $h_n$ commute with the degeneracy maps $s_n^i$ for
$0 \leq i \leq n$ and 
with the face maps $d_n^i$ for $0 \leq i<n$. For $d_n^n$ we get that:
$$f_{n-1} \circ d_n^n (a,v_1,v_2, \dots, v_n) = 
(v_{n-1}^{-1}\cdots v_2^{-1}v_1^{-1}av_1,v_2, \dots, v_{n-1})$$
$$d^n_n \circ f_n (a,v_1,v_2, \dots, v_n) =
(v_{n-1}^{-1}\cdots v_2^{-1}v_1^{-1}av_1,v_2, \dots, v_{n-1}).$$
\end{proof}

Moreover, both simplicial spaces  $\wedge \Gg_*$ and $\YY_*(\Gg)$ can be given
a structure of cyclic space (see \ref{definitioncyclic}) in such a way that they remain isomorphic as cyclic spaces.

\begin{rem}
The space $\YY_*(\Gg)$ is also known by the name of {\it cyclic nerve} of the groupoid $\Gg$.
\end{rem}

The cyclic structure associated to $\wedge \Gg_*$ is:
\begin{eqnarray*}
t_n : \wedge \Gg_n & \to & \wedge \Gg_n \\
(a, v_1, \dots, v_n) & \mapsto & ( (v_n^{-1} \cdots v_1^{-1} a v_1 \dots v_n), 
(v_n^{-1} \cdots v_1^{-1} a), v_1, \dots, v_{n-1})
\end{eqnarray*}

$$
\xymatrix{
         \circ \ar@(ul,ur)[]^a \ar@/^/[r]^{v_1} & \circ \ar@{.}[r] & \circ \ar@/^/[r]^{v_{n-1}}
& \circ \ar@/^/[r]^{v_n} & \circ}
\mapsto 
\xymatrix{
         \circ  \ar@/^/[r]^{v_1} & \circ \ar@{.}[r] & \circ \ar@/^/[r]^{v_{n-1}} &
 \circ & \circ \ar@(ul,ur)[]^{v_n^{-1} \cdots v_1^{-1}a v_1\cdots v_n} 
\ar@/^1pc/[llll]^{v_n^{-1} \cdots v_1^{-1}a} }
$$

and the cyclic structure of $\YY_*(\Gg)$ is:
\begin{eqnarray*}
t_n : \YY_n(\Gg) & \to & \YY_n(\Gg) \\
(g_0, g_1, \dots, g_n) & \mapsto & (g_n, g_0, g_1, \dots, g_{n-1})
\end{eqnarray*}
It is easy to verify that the simplicial maps $f_*$ and $h_*$ become cyclic maps,  hence
\begin{proposition} \label{isoofcyclic}
The maps $f_* : (\wedge \Gg_*, t_*) \to (\YY_*(\Gg),t_*)$ and
$h_* : (\YY_* (\Gg), t_*) \to (\wedge \Gg_*,t_*)$ are isomorphisms of cyclic spaces.
\end{proposition}

The cyclic structure on a simplicial space $\XX_*$ allows one to define a circle action
on the geometrical realization $|\XX_*|$ as is shown in \cite{BurgheleaFiedorowicz}.
The construction is explained in the Appendix, Prop. \ref{circleactioncyclic}.

Going back to the configuration space $\CC_{S^1}(\Gg)$ we have

\begin{proposition} \label{conf=inertia}
The configuration space $\CC_{S^1}(\Gg)$ is homeomorphic to the geometrical
realization $|\YY_*(\Gg)|$. Moreover, the homeomorphism is $S^1$-equivariant.
\end{proposition}

\begin{proof}

The argument is very similar to the one in Proposition \ref{propositionBG}. We 
define a bijective map $\varphi : \CC_{S^1}(\Gg) \to |\YY_*(\Gg)|$ using the fact that as sets
$$|\YY_*(\Gg)| = \bigsqcup (\YY_n(\Gg) - {\rm Deg} \ \YY_n(\Gg)) \times {\rm Int} \ \Delta_n.$$

For a configuration $\alpha \in \CC_{S^1}(\Gg)$ with $${\rm supp}(\alpha)=\{x_0, x_1, \dots, x_n\}, \ \
0 \leq x_0 < x_1 < \cdots < x_n < 1$$ and $g_i:=\alpha(x_i)$ we will associate the point 
\begin{eqnarray*} \label{def-phi}
\varphi(\alpha) = \left\{
\begin{array}{ccc}
(\ident \source (g_0), g_0, g_1, \dots, g_n ; x_0, x_1-x_0, \dots, 1-x_n) & {\rm if} & x_0 \neq 0 \\
 ( g_0, g_1, \dots, g_n ;  x_1, x_2-x_1 \dots, 1-x_n) & {\rm if} & x_0 = 0.
\end{array} \right.
\end{eqnarray*}
This function is clearly a homeomorphism.

To prove that $\varphi$ is $S^1$-equivariant we will mimic to proof given in \cite{Burghelea1}. Let's recall
the rotation map $\rho : S^1 \times \CC_{S^1}(\Gg) \to \CC_{S^1}(\Gg)$ defined in \ref{rotationmap}
and the circle action $\mu : S^1 \times |\YY_*(\Gg)| \to |\YY_*(\Gg)|$ explained in \ref{circleactioncyclic}; we
want to verify the commutativity of the following diagram

\begin{eqnarray*}
\xymatrix{
        &  S^1 \times \CC^{S^1}(\Gg)  \ar[r]^\rho  &   \CC_{S^1}(\Gg) & \\
	 |S^1_*| \times |\YY_*(\Gg)| \ar[ur]^{i \times \varphi^{-1}} & |S^1_* \times \YY_*(\Gg)| \ar[l]_j
    & |S^1_* \stackrel{\sim}{\times} \YY_*(\Gg)| \ar[l]_{h} \ar[r]^{|ev|} & |\YY_*(\Gg)|. \ar[ul]_{\varphi^{-1}}}
\end{eqnarray*}
We take a point in $|S^1_* \stackrel{\sim}{\times} \YY_*(\Gg)|$ written as $((\tau^r_n,x);(u_0,\dots,u_n))$ 
$0 < u_i \leq 1$, $\sum u_i = 1$ and $x=(g_0,\dots,g_n)$ with $g_s \notin \ident(\Gg_0)$ for $s= 1,\dots,n$.

In order to make the following argument more understandable, let us recall what the configuration $\alpha$
associated to the pair $(x;u_0, \dots, u_n) \in \YY_*(\Gg)$ given by the map $\varphi^{-1}$ is.
As the $u_i$'s are positive then $\alpha(0)=g_0$, $\alpha(u_0)= g_1$, \dots, $\alpha(u_0 + \cdots + u_{n-1})=g_n$.

Going back to the diagram we have following set of equalities
\begin{eqnarray*}
\lefteqn{(i\times \varphi^{-1})\cdot j \cdot h ((\tau^r_n,x);(u_0,\dots,u_n))}\\
 & = & (i\times \varphi^{-1})\cdot j ((\tau^{n-r+1}_n,x);(u_r, u_{r+1}, \dots,u_n, u_0,\dots,u_{r-1}))\\
&=& (i \times \varphi^{-1})[(\tau^1_1;s,t),(x;u_r, u_{r+1}, \dots,u_n, u_0,\dots,u_{r-1})]\\
&=& \left( i(\tau^1_1;s,t);\varphi^{-1}(x;u_r, u_{r+1}, \dots,u_n, u_0,\dots,u_{r-1}) \right)\\
&=& (s,\gamma)
\end{eqnarray*}
with $s=u_r + u_{r+1} + \cdots + u_n$ , $t=u_0 +\cdots + u_{r-1}$, and $\gamma$ the configuration given by
$\gamma(0)=g_0$, $ \gamma(u_r) =g_1$, $\gamma(u_r + u_{r+1}) = g_2$,  $ \gamma(u_0+\cdots +u_r) = g_{n-r+1}$,
 $ \gamma(u_r+ \cdots + u_n+u_0+\cdots + u_{r-2})  = g_n$.

Then the configuration $\lambda:= \rho(s,\gamma)$ is given by $\lambda(u) = \gamma(s+u)$; therefore
 $\lambda(0) =\gamma(s)=g_{n-r+1}$,
$\lambda(u_0 +\cdots + u_{r-1}) = \gamma(s+t) = \gamma(0) = g_0$, $\lambda(u_0 +\cdots+u_{n-1}) =
 \gamma(u_r+\cdots +u_{n-1}) = g_{n-r}$.

On the other side of the diagram we have
\begin{eqnarray*}
\lefteqn{ \varphi^{-1} \cdot |ev| ((\tau^r_n,x);(u_0,\dots,u_n))}\\
& =& \varphi^{-1}(t^r_n x; u_0,\dots,u_n) \\
& = & \varphi^{-1}((g_{n-r+1},g_{n-r+2}, \dots ,g_n, g_0, \dots, g_{n-r}); u_0,\dots,u_n)\\
& = & \beta  
\end{eqnarray*}
with $\beta(0) = g_{n-r+1}$, $\beta(u_0)=g_{n-r+2}$, $\beta(u_0 +\cdots + u_{r-1})=g_0$, $\beta(u_0+\cdots + u_{n-1}) = g_{n-r}$.
Hence $\beta = \lambda$ and the diagram is commutative.
\end{proof}

We can conclude
\begin{cor}
The classifying space of the inertia groupoid $B \wedge\Gg = |\YY_*(\Gg)|$ is $S^1$-homeomorphic
to the configuration space $\CC_{S^1}(\Gg)$.
\end{cor}

\section{The ghost loop space} \label{sectionghostloopspace}

We will postpone formal definitions for a while and we will start instead by explaining what a  ``ghost loop'' is.

An object in the inertia groupoid can be seen as a morphism of groupoids
$$ \Psi:\bar{\integer} \to \Gg$$
where we think of $\bar{\integer}$ as the groupoid with one object $\{ *\}$ and 
whose morphisms are the integers (see \cite[Corollary 3.6.4]{LupercioUribe2}).
If we get the geometrical realization of the morphism $\Psi$, 
$B\Psi : B\integer  \to B \Gg$ we obtain a map $\psi: S^1 \to B \Gg$. Now, is any map
$S^1 \to B\Gg$ homotopic to the realization of a morphism $\Psi$? The answer of course
depends on the groupoid $\Gg$, but in general this is not the case. This section
is devoted to give a very down to earth, geometrical description of the realization
of such morphisms.

The loop $\psi = B\Psi$ will be \emph{anchored} at the point $x:=\Psi(*) \in \Gg_0$ and its
image will live on the fiber $\pi^{-1}([x]) \simeq BG_x$ where
$$\pi: B\Gg \to \coarse{\Gg}$$
is the projection of $B \Gg$ in the coarse moduli space $\coarse{\Gg}:=\Gg_0/\Gg_1$ (where $\Gg_1$ provides
the equivalence relation for the elements in $\Gg_0$), $[x]$ is the equivalence class of $x$ in
$\coarse{\Gg}$ and 
$$G_x := \{g \in \Gg_1 \ | \  \source(g)=\target(g) = x\}$$ 
is the isotropy group at $x$.
\begin{center}
\includegraphics[height=3.5in]{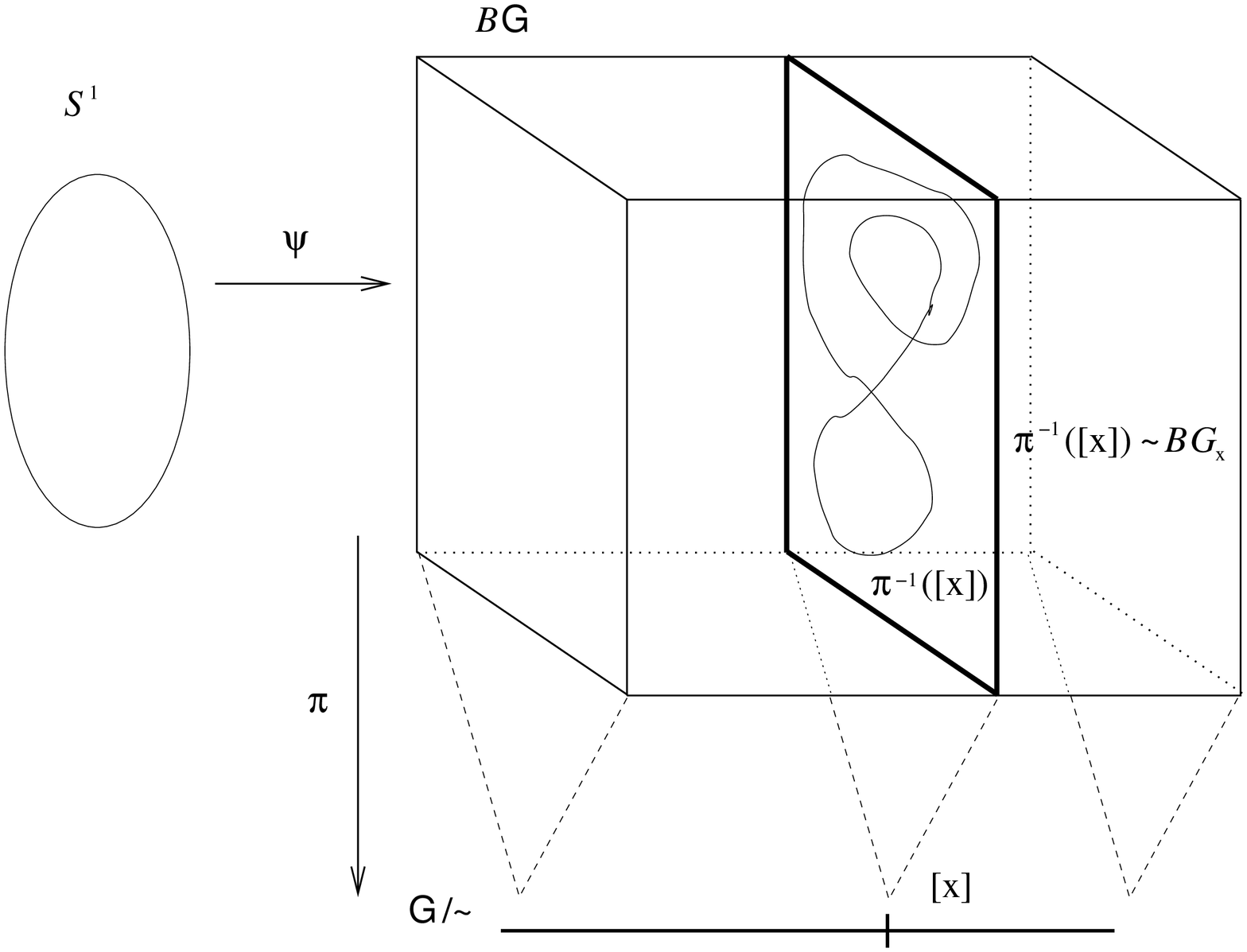}

\end{center}

These are the loops that we will call ``ghost loops'' (our terminology is motivated by string theory.)
After this brief introduction we are ready to define the spaces we need.

Let $\LL B \Gg:= {\rm Map}(S^1,B\Gg)$ be the space of loops in $B\Gg$ with the natural
compact-open topology.

\begin{definition}
The space $\smallloop B\Gg$ of ghost loops of $B\Gg$ is the set of maps $f \in  \LL B \Gg$
such that the image of $f$ lies in one fiber of $\pi$, that is, if it exists $x \in \Gg_0$ with
$${\rm Im}\ f \subset \pi^{-1}([x])$$
and $[x] \in \coarse{\Gg}$.
\end{definition} 
It is equivalent to say that the following diagram commutes:
$$
\xymatrix{
         S^1 \ar[r]^{f} \ar[d] & B\Gg \ar[d]^{\pi}\\
\{[x]\} \ar@{^{(}->}[r] & \coarse{\Gg} }$$

From the isomorphism of Proposition \ref{propositionBG} and the evaluation map defined in formula
(\ref{evaluationmap}) we obtain a map $\mu : \CC_{S^1}(\Gg) \times S^1 \to B\Gg$. This map induces
 a natural $S^1$ equivariant map 
$$\nu : \CC_{S^1}(\Gg) \to \LL B \Gg$$ 
where the action on the l.h.s. is explained in formula (\ref{rotationmap}) and the one in the r.h.s is the 
obvious one. But the image of $\nu $ in $\LL B \Gg$ lies  in the ghost loops space.
\begin{lemma}
${\rm Im}(\nu) \subset \smallloop{B\Gg}$.
\end{lemma}

\begin{proof}
Let $\alpha$ be a configuration in $\CC_{S^1}(\Gg)$ with $y \in {\rm supp}(\alpha)$ and $g=\alpha(y)$.
If $x = \source(g)$ it is clear that $\pi ( \mu(\alpha, \theta)) = [x]$ for any $\theta \in S^1$. Hence
${\rm Im}\ \nu(\alpha) \subset \pi^{-1}([x])$.
Without confusion we can assume that 
$$\nu: \CC_{S^1}(\Gg) \to \smallloop{B\Gg}.$$
\end{proof}

Let's consider the following commutative diagram
$$\xymatrix{
\CC_{S^1}(\Gg) \ar[r]^{\varepsilon_0} \ar[d]^\nu & \CC_R(\Gg) \ar[d]^{\cong}\\
\smallloop{B \Gg} \ar[r]^{ev_0} & B\Gg}
$$
where $\varepsilon_0$ is the restriction of the map $\varepsilon $ to $\CC_{S^1}(\Gg) \times \{0\}$,
$ev_0$ is the evaluation of a loop at $0$, i.e. for $f \in \LL B \Gg$, $ev_0(f):=f(0) \in B\Gg$ and
$\cong$ is the $S^1$-equivariant homeomorphism of Proposition \ref{propositionBG}.

\begin{theorem}
 \label{homotopyonfibers}
The map $\nu: \CC_{S^1}(\Gg) \to \smallloop{B\Gg}$ is $S^1$-equivariant and induces a weak homotopy equivalence
on the fibers $$\nu_\alpha : {\varepsilon_0}^{-1}(\alpha) \stackrel{\simeq}{\longrightarrow} {ev_0}^{-1}(\alpha)$$
for every $\alpha \in \CC_R(\Gg)$ . 
\end{theorem}

\begin{proof} Let $\alpha$ be as in definition \ref{definitionconfiguration}; a map $\alpha:R=(0,1) \to \Gg_1$ with support
${\rm supp} (\alpha) := \{t_1, t_2, \dots, t_n \}$, $g_i : = \alpha(t_i)$ and $\target(g_i)=\source(g_{i+1})$. The space
${\varepsilon_0}^{-1}(\alpha)$ will consist of all the configurations on the circle $\beta \in  \CC_{S^1}(\Gg)$ such that when
restricted to the interval $(0,1)$ matches $\alpha$, i.e $\beta|_{(0,1)} = \alpha$. If we call $x:=\source(g_1)$ and

 $y:=\target(g_n)$
it is clear that all the possible values of $\beta(0)$ lie in the set
$${\rm Hom}_\Gg(x,y):=\{g \in \Gg_1 | \source(g)=y \ \& \ \target(g)=x \},$$
and is easy to see that ${\varepsilon_0}^{-1}(\alpha)\cong {\rm Hom}_\Gg(x,y)$.

There is a canonical homeomorphism between the space ${\rm Hom}_\Gg (x,y)$ and the isotropy group $G_x$ at $x$ given
by the following map
\begin{eqnarray*}
{\rm Hom}_\Gg (x,y) & \stackrel{\simeq}{\longrightarrow} & G_x \\
g & \mapsto & g_1g_2 \dots g_ng
\end{eqnarray*}

Recall now that  the space ${ev_0}^{-1}(\alpha)$ consist of all loops $f : S^1 \to B \Gg$ in  $\smallloop{B \Gg}$ such that
$f(0) = \alpha$. We have already seen that the image of these loops lie in $\pi^{-1}([x])$ which is homotopic to
$B G_x$.
Let $\Gg^x_0$ be the orbit of $x$ under the groupoid $\Gg$ and let's construct the groupoid $\Gg^x$ whose
morphisms are given by the set $\Gg^x_1$ which makes the following diagram into a cartesian square
$$\xymatrix{
\Gg^x_1 \ar[r]  \ar[d]^{\source \times \target} & \Gg_1 \ar[d]^{\source \times \target}\\
\Gg^x_0 \times \Gg^x_0 \ar[r] & \Gg_0 \times \Gg_0}$$
We can identify $\Gg^x$ as a subgroupoid of $\Gg$ and therefore $B\Gg^x = \pi^{-1}([x])$.
Then ${ev_0}^{-1}(\alpha)$ is the space of loops in $B\Gg^x$ based at $\alpha$ (i.e. $\Omega B \Gg^x$).

So we obtain the following maps
$$G_x  \stackrel{\cong}{\longrightarrow} {\rm Hom}_\Gg (x,y) \stackrel{\cong}{\longrightarrow} {\varepsilon_0}^{-1}(\alpha)
\stackrel{\nu_\alpha}{\longrightarrow} {ev_0}^{-1}(\alpha) \stackrel{\cong}{\longrightarrow} \Omega B \Gg^x
\stackrel{\simeq}{\longrightarrow} \Omega B G_x$$
where the composition of all of them $G_x \to \Omega BG_x$ is the weak homotopy equivalence induced by
by the embedding $\{ h \} \times \Delta_1 \subset | \YY_*(\Gg)| \cong \CC _{S^1}(G_x) \simeq \LL BG_x$ for $h \in G_x$.

The fact that $\nu$ is $S^1$-equivariant follows from its definition.
\end{proof}

We conjecture that the previous information is sufficient in order for $\nu$ to induce a weak homotopy equivalence.
In what follows we will prove that this is the case for orbifolds. The fact that orbifolds
are locally quotients of manifolds by finite groups allows us to do so. Let us briefly recall some known
properties of orbifolds that can be found in \cite{Moerdijk2002, MoerdijkPronk}.

\subsection{Orbifolds}

For every orbifold $(X,\UU)$ it can be chosen an atlas $\UU=(\pi_i, U_i, G_i, V_i)_i$ with the following properties
 \cite[Cor. 1.25]{MoerdijkPronk}:
\begin{itemize}
\item For every chart both $U_i$ and $V_i$ are contractible,
\item The intersection of finitely many charts in $\UU$ is either empty or again a chart in $\UU$,  and
\item The coordinates can be chosen  so that each $U_i$ is an euclidean ball and the $G_i$'s act linearly.
\end{itemize}

The groupoid $\Gg$ associated to such orbifold atlas (as shown in \cite{LupercioUribe1, MoerdijkPronk}) is {\it Leray}
 i.e. the spaces $\Gg_n$, $n\geq 0$ are diffeomorphic to a disjoint union of contractible open sets. In this
 case the coarse moduli space $\coarse{\Gg}$ is equal to $X$ and the restriction $\Gg|_{V_i}$ of the groupoid
  $\Gg$ to $V_i$ via the quotient map $\Gg \to X$, is Morita
 equivalent to the groupoid $\Global{U_i}{G_i}$ (the groupoid constructed via the action of $G_i$ on $U_i$).
In the case that $U$ is an euclidean ball and the action of $G$ is linear we have the next result.

\begin{lemma}
The map $\nu : \CC_{S^1}(\Global{U}{G}) \to \smallloop B \Global{U}{G}$ induces a weak homotopy equivalence.
\end{lemma}
\begin{proof}
The key fact in what follows is that the ball can be equivariantly contracted to the origin.

Let $h_t:U \to U$ be the homotopy $u \mapsto ut$ that contracts the ball to the origin and let $h_t'$ be the homotopy
induced  in the ghost loop spaces
$$h_t' : \smallloop B \Global{U}{G} \To \smallloop B \Global{U}{G}.$$
Then $h_1'$ is the identity and $h_0'$ factors through $\smallloop BG= \LL BG$; then it induces a weak
homotopy equivalence $\smallloop B \Global{U}{G} \stackrel{\simeq}{\To} \LL BG$.

As the action of $G$ is linear, the same procedure could be applied to the inertia groupoid $\wedge \Global{U}{G}$. We
obtain the following commutative diagram

$$
\xymatrix{
          |\wedge \Global{U}{G}| \ar[r]^{\cong} \ar[d]^{\simeq} & \CC_{S^1}(\Global{U}{G}) \ar[d]^{\simeq} \ar[r]^\nu &
	  \smallloop B \Global{U}{G} \ar[d]^\simeq \\
          |\wedge G | \ar[r]^\cong & \CC_{S^1}(G) \ar[r]^{\simeq}_{\nu_G}& \LL B G }$$
The map $\nu_G$ induces a weak homotopy equivalence as was described in \ref{homotopyonfibers}; the result follows.
\end{proof}

Using the orbifold atlas $\UU$ we can associate to $X$ a simplicial space $X_*$ constructed from the open cover
 $\{V_i\}_{i \in I}$ as follows. The $X_{n}$'s will consist of the disjoint union of the intersections of $n+1$ sets from
 the open cover and the face maps $d_j :X_{n+1} \to X_{n}$ will be the natural inclusions, in other words,
 if we call $V_{i_0 \dots i_n} := V_{i_0} \cap \cdots \cap V_{i_n}$ then
 $$X_n := \bigsqcup_{(i_0,\dots ,i_n) \in I^{n+1}} V_{i_0\dots i_n}$$
 and
 $$ d_j :V_{i_0\dots i_n} \hookrightarrow V_{i_0\dots \widehat{i_j} \dots i_n}.$$

It is a result due to Segal that $|X_*|$, the geometrical realization of $X_*$, is weakly homotopy equivalent to $X$.

Now let's define the maps
$\widetilde{\varepsilon_0}= \pi \circ \varepsilon_0 : \CC^{S^1}(\Gg) \to X$ and
$\widetilde{ev_0}= \pi \circ ev_0 : \smallloop B\Gg \to X$,
where $\pi : B \Gg \to \coarse{\Gg} =X$ is the natural projection. Let $Y_*$ and $Z_*$ be the simplicial spaces associated to $X_*$
and the maps $\widetilde{\varepsilon_0}$ and $\widetilde{ev_0}$ respectively. This means
$$Y_n:=  \widetilde{\varepsilon_0}^{-1} (X_{n})  \ \ \ \ \mbox{and} \ \ \ \ \ Z_n:=  \widetilde{ev_0}^{-1}  (X_{n})$$
with the natural inclusions as face maps.

By the previous lemma $\nu$ induces a weak homotopy equivalence $\nu_n : Y_n \stackrel{\simeq}{\to} Z_n$
and at the same time a map of simplicial spaces $\nu_* :Y_* \to Z_*$. Then taking geometric realizations
$$\CC_{S^1}(\Gg) \simeq |Y_*| \stackrel{\simeq}{\To} |Z_*| \simeq \smallloop B \Gg$$
we obtain a weak homotopy equivalence\footnote{The simplicial sets $Y_*, Z_*$ are constructed from the open covers
$\{\widetilde{\varepsilon_0}^{-1}(V_i)\}_{i\in I},\{\widetilde{ev_0}^{-1}(V_i)\}_{i\in I} $ of $\CC_{S^1}(\Gg)$ and
$\smallloop B \Gg$ respectively.}
 between the configuration space and the ghost loop space. So we get

\begin{theorem} \label{w.h.e.conf=ghost}
For $\Gg$ an orbifold groupoid the map $\nu: \CC_{S^1}(\Gg) \to \smallloop B \Gg$ induces a weak homotopy
equivalence.
\end{theorem}

It is easy to see that all the previous procedure can be carried out after performing the $S^1$
Borel construction on the configuration space and on the ghost loop space. This because the maps $\nu, \nu_n$ and the
ones on the previous lemma are all $S^1$ equivariant. Then

\begin{cor}
The map $\nu$ induces a weak homotopy equivalence between $ES^1 \times_{S^1} \CC_{S^1}(\GG)$ and
$ES^1 \times_{S^1} \smallloop B \Gg$.
\end{cor}

So we have the following isomorphisms in homology

\begin{cor} \label{isoinhomology}
$$H_*(\CC_{S^1}(\Gg), \integer) \stackrel{\cong}{\longrightarrow} H_*(\smallloop B\Gg, \integer)$$
$$H^{S^1}_*(\CC_{S^1}(\Gg), \integer) \stackrel{\cong}{\longrightarrow} H^{S^1}_*(\smallloop B\Gg, \integer).$$
\end{cor}
Now we will relate the previous study to other homology theories.

\section{Hochschild and cyclic homology}
\subsection{Preliminaries}
\begin{definition}
A cyclic $R$-chain complex $A$ is a functor from $\Lambda^{op}$ to $R$-chain complexes (see Remark \ref{cycliccategory}); it
is determined by the chain complexes $A_n:=A(\mathsf{n})$ and the morphisms $d_n^i: A_n \to A_{n-1}$, $S^i_n: A_n \to A_{n+1}$
 and
$t_n : A_n \to A_n$. One defines the {\it Hochschild homology} $HH_*(A_\bullet)$
as the total homology of the chain complex of chain complexes
\begin{eqnarray} \label{HHcomplex}
\xymatrix{
    A_0 & A_1 \ar[l]_b & A_2 \ar[l]_b & \cdots \ar[l]_b & A_n \ar[l]_b \ar[l]_b & \cdots \ar[l]_b  }
\end{eqnarray}
and the {\it cyclic homology} $HC_*(A_\bullet)$ as the homology of the chain bicomplex of chain of complexes
\begin{eqnarray} \label{HCcomplex}
\xymatrix{
    \vdots \ar[d]_b & \vdots \ar[d]_{-b'} & \vdots \ar[d]_b & \vdots \ar[d]_{-b'} & \\
    A_2 \ar[d]_b & A_2 \ar[d]_{-b'} \ar[l]_{1-T} & A_2 \ar[d]_b \ar[l]_{N}  & A_2 \ar[l]_{1-T} \ar[d]_{-b'} & \cdots  \ar[l] \\
    A_1 \ar[d]_b & A_1 \ar[d]_{-b'} \ar[l]_{1-T} & A_1 \ar[d]_b \ar[l]_{N}  & A_1 \ar[l]_{1-T} \ar[d]_{-b'} & \cdots  \ar[l] \\
    A_0 & A_0 \ar[l]_{1-T} & A_0 \ar[l]_{N}  & A_0 \ar[l]_{1-T} & \cdots  \ar[l]}
\end{eqnarray}
with $b, b': A_n \to A_{n-1}$, $T, N:A_n \to A_n$ given by
\begin{eqnarray*}
b= \sum_{i=0}^n (-1)^{i} d^i_n, \ \ \  b'= \sum_{i=0}^{n-1} (-1)^{i} d^i_n, \ \ \ T=(-1)^n t_n, \ \ \ N= 1+ T + \cdots + T^n.
\end{eqnarray*}
\end{definition}
As the complex
\begin{eqnarray*}
\xymatrix{
    A_0 & A_1 \ar[l]_{-b'} & A_2 \ar[l]_{-b'} & \cdots \ar[l]_{-b'} & A_n \ar[l]_{-b'} \ar[l]_{-b'} & \cdots \ar[l]_{-b'}  }
\end{eqnarray*}
is acyclic one obtains an exact sequence
\begin{eqnarray*}
\xymatrix{
   \cdots \ar[r] & HH_*(A_\bullet) \ar[r]^I & HC_*(A_\bullet) \ar[r]^S & HC_{*-2}(A_\bullet) \ar[r] & HH_{*-1}(A_\bullet)
   \ar[r] & \cdots}
\end{eqnarray*}
with $I$ induced by the inclusion of \ref{HHcomplex} in \ref{HCcomplex} as the first column and $S$ induced by the
canonical projection of \ref{HCcomplex} on itself which sends the first two columns to zero and the $(k+2)$-th column
to the $k$-th one. This sequence is called the {\it Connes sequence} of $A_\bullet$.

\begin{definition} \label{defperiodic}
The homology of the inverse limit of complexes induced by $S$ is denoted $HP_*(A_\bullet)$ and
 is called {\it periodic homology}; it
is a $\integer / 2 \integer$ graded theory.
\end{definition}

\begin{example}
If $(X_*,t_*)$ is a cyclic space, we consider $\XX$ to be the composition of $(X_*,t_*) : \Lambda^{op} \to Top$ with
$C_*(\cdot; R) : Top \to R-chain \ complexes$, where $C_*(X;R)$ is the normalized singular chain complex of $X$ with
coefficients in $R$.
One denotes by $HH_*((X_*,t_*);R)$ resp.  $HC_*((X_*,t_*);R)$ the Hochschild homology $HH_*(\XX_\bullet)$ resp. cyclic homology
$HC_*(\XX_\bullet)$.
\end{example}
Is a result  in \cite{BurgheleaFiedorowicz} that for $R$ a commutative ring with unit
\begin{proposition} \label{resultBurgheleaFiedorowicz}
For any cyclic space $(X_*,t_*)$ there exist natural isomorphisms
\begin{eqnarray*}
HH_*((X_*,t_*);R) & \stackrel{\cong}{\To} & H_*(\Vert X_*\Vert ; R) \\
HC_*((X_*,t_*);R) & \stackrel{\cong}{\To} & H_*(\Vert X_*,t_* \Vert ; R)
\end{eqnarray*}
which identifies Connes' homomorphism $S:HC_*((X_*,t_*);R) \to HC_{*-2}((X_*,t_*);R)$ with the Gysin homomorphism
$S:H_*(\Vert X_*,t_* \Vert;R) \to H_{*-2}(\Vert X_*,t_*\Vert ;R)$ of the fibration
$\Vert X_* \Vert \to \Vert X_*,t_* \Vert \to BS^1$.
\end{proposition}

\subsection{Orbifolds}

Applying the previous formalism to orbilfod groupoids we obtain the following results,

\begin{theorem}
Let $\Gg$ be an orbifold groupoid and $\wedge \Gg$ its inertia groupoid. Then there are
canonical isomorphisms
\begin{eqnarray*}
HH_*((\wedge \Gg_*,t_*)) \stackrel{\cong}{\To} H_*(|\YY_*(\Gg)|) & \stackrel{\cong}{\To} & H_*(\CC_{S^1}(\Gg))
\stackrel{\cong}{\To} H_*(\smallloop B \Gg)\\
  HC_*((\wedge \Gg_*,t_*)) \stackrel{\cong}{\To} H^{S^1}_*(|\YY_*(\Gg)|) & \stackrel{\cong}{\To} & H^{S^1}_*(\CC_{S^1}(\Gg))
\stackrel{\cong}{\To} H^{S^1}_*(\smallloop B \Gg)
\end{eqnarray*}
\end{theorem}

\begin{proof}
The LHS isomorphisms follow from  \ref{resultBurgheleaFiedorowicz} and 
\ref{fibration}, the isomorphisms in the middle from \ref{conf=inertia} and the RHS ones from \ref{isoinhomology}
\end{proof}

\subsection{Chen-Ruan Orbifold cohomology} Motivated by the theory of quantum cohomology Chen
and Ruan \cite{ChenRuan} defined a remarkable cohomology
ring associated to a smooth complex orbifold. One of the peculiarities of the theory is that in certain interesting cases
is isomorphic to the cohomology ring of the smooth space which resolves the singularities of the
orbifold \cite{Ruan, Uribe, Fantechi}, namely a crepant resolution.

As we mentioned in the introduction, orbifolds can be seen as a special kind of topological groupoids
 i.e. smooth, stable and \'etale \cite{LupercioUribe1,Moerdijk2002, MoerdijkPronk}.
The twisted sectors $(X_{(g)}, g \in T_1)$ of an orbifold $X$ \cite{ChenRuan,Satake}, on which the definition
 of the orbifold cohomology relies,
are known to be equivalent
to the inertia groupoid of the orbifold (seen as a groupoid) \cite{AbramovichGraberVistoli,LupercioUribe2, Moerdijk2002}.
 To each
connected component of the twisted sectors there is associated a rational number $\iota_{(g)}$ 
which is called {\it degree shifting number} \cite{ChenRuan} 
(also known as {\it age} \cite{Reid} ), and the orbifold cohomology is defined as the direct sum of the 
cohomology of the twisted sectors whose degree is shifted up by twice the degree shifting number, i.e.
$$H_{\mathrm{orb}}^*(X,\complex) \stackrel{def}{=} \bigoplus_{(g)} H^{*-2\iota_{(g)}}(X_{(g)}, \complex).$$
 For $SL$-orbifolds (e.g. Gorenstein orbifolds,) namely orbifolds for which the local action at every point
  can be seen as an element of $SL(d,\complex)$, the shifting number is an integer. In the case of compact
   orbifolds , i.e proper \'{e}tale groupoids whose coarse moduli space is compact, we get the following
    immediate consequence of the results of the previous section.

\begin{cor}
For $\Gg$ a compact complex $SL$-orbifold the groups
$$\prod_{m \in n + 2 \integer} HH_m((\wedge \Gg_*, t_*) ; \complex) 
\cong \prod_{m \in n + 2 \integer}H_{\mathrm{orb}}^m(\Gg ; \complex)$$
are isomorphic.
\end{cor}
\begin{proof}
Here we use the fact that the homological dimensions of the connected components of the inertia
 groupoid are all even (because the orbifold is
complex) as well as the shiftings on the degree.
\end{proof}

Finally using the Chern character for orbifold $K$-theory \cite{AdemRuan, Moerdijk2002} and reduced orbifolds (i.e.
for every $x \in \Gg_0$ the induced homomorphism $G_x \to Aut(\Gg_0)$ is injective) ,we can obtain the following isomorphism.

\begin{theorem}
Let $\Xx$ be a compact complex reduced $SL$-orbifold and $\Gg$ a groupoid associated to it. Then
 the periodic homology of $\wedge \Gg$ is
isomorphic to the Chen-Ruan orbifold cohomology of $\Xx$ (but $\integer/2\integer$-graded)
$$K^n_{\mathrm{orb}}(\Xx) \otimes \complex  \cong \prod_{m \in n + 2 \integer} H_{\mathrm{orb}}^m(\Xx, \complex)\cong 
HP_n((\wedge \Gg_*,t_*), \complex).$$
\end{theorem}

\begin{proof}

The first isomorphism is given by the Chern character map of orbifold $K$-theory \cite{AdemRuan, Moerdijk2002},
 and the second follows from
the previous corollary and the fact that for orbifolds there is an isomorphism (Cor. \ref{periodic=hochschild})
$$HP_n((\wedge \Gg_*, t_*), \complex) \cong \prod_{m \in n + 2 \integer} HH_m((\wedge \Gg_*,t_*), \complex).$$

\end{proof}

\begin{proposition}
Let $\Xx$ be an orbifold  and $\Gg$ a groupoid associated to it. Then the Connes exact sequence split into the 
short exact sequences
\begin{eqnarray*}
0 \to HH_k((\wedge \Gg_*, t_*); \rational) \to HC_k((\wedge \Gg_*, t_*); \rational) 
\to HC_{k-2}((\wedge \Gg_*, t_*);\rational) \to 0
\end{eqnarray*}  
for every $k \in \integer$.
\end{proposition}
\begin{proof}
In rational homology the action of the circle $S^1$ on $|\wedge \Gg_*|$ is trivial. This because the 
isotropy is finite and the classifying
spaces of finite groups are rationally trivial. Then the spectral sequence associated to the fibration
 of theorem \ref{fibration} reduces to
the isomorphism
$$H^{S^1}_*(|\wedge \Gg_*|; \rational) \cong H_*(|\wedge \Gg_*| ; \rational ) \otimes H_*^{S^1}(*; \rational).$$
Applying theorem \ref{resultBurgheleaFiedorowicz} the result follows. 
\end{proof}
In the case of an orbifold the periodic homology is the inverse limit of the maps $S: HC_n \to HC_{n-2}$,
 therefore we can easily see that
\begin{cor} \label{periodic=hochschild}
For $\Xx$ a $\mathrm{d}$-dimensional orbifold, $\Gg$ a groupoid associated to it, $k=0,1$ and any $l$ 
such that $ 2l \geq d$ then
$$HP_k((\wedge \Gg_*,t_*);\rational) \cong HC_{2l+k}((\wedge \Gg_*,t_*);\rational) \cong 
\prod_{n=0}^{|\frac{d}{2}|} HH_{2n +k}((\wedge \Gg_*,t_*);\rational).$$
\end{cor}

\begin{rem}
The previous isomorphisms are additive. We will return to the study of the multiplicative structures in a future paper.
\end{rem}
 
\section{Appendix}

This section will be devoted to give a summary of some of the properties of simplicial
and cyclic spaces. We will be mostly based on the papers \cite{BurgheleaFiedorowicz, FiedorowiczLoday} 
but we recommend the reader the book of Connes \cite{Connes} for a friendly
description of the subject.

\subsection{Simplicial spaces}

Let $Ord$ be the category whose objects are the sets $\mathsf{n}=\{0,1,\dots,n\}$ for $n \in \naturals$
and whose morphisms are the order preserving maps $\mathsf{n} \to \mathsf{m}$. The face maps $\partial^n_i 
:\mathsf{n} -1 \to \mathsf{n}$ which skip $i$, and the degeneracy maps  $\lambda^n_i 
:\mathsf{n}+1 \to \mathsf{n}$ which repeat $i$, generate all the morphism of the category.

\begin{definition}
A {\it simplicial space} $X$ is a contravariant functor from $Ord$ to the
category of topological spaces $Top$.
\end{definition}

It can also be described as $(X_*,d^i_n,s^i_n)$ where $X_n:=X(\mathsf{n})$, $d^i_n:=X(\partial^n_i)$,
$s^i_n:=X(\lambda^n_i)$. Hence $d_n^i:X_{n} \to X_{n-1}$ and $s^i_n:X_{n} \to X_{n+1}$ with
\begin{eqnarray*}
d^i_{n-1} d^j_n & =  & d^{j-1}_{n-1} d^i_n \  \ \ \ \ \ \ \ \ \  \ \ \ \ i<j\\
s^i_{n+1} s^j_n & =  & s^{j+1}_{n+1} s^i_n \ \ \ \ \ \  \ \ \ \ \ \ \ \  i \leq j\\
d^i_{n+1} s^j_n & =  & \left\{
\begin{array}{cc}
s^{j-1}_{n-1} d^i_n &  \ \ \ \ \ \ i<j\\
1 & \ \ \ \ \ \ i=j,j+1\\
s^j_{n-1} d^{i-1}_n & \ \ \ \ \ \ i>j+1
\end{array} \right.
\end{eqnarray*}

A simplicial map between  simplicial spaces $f_*:X_* \to Y_*$ is a family of continuous maps
$f_n :X_n \to Y_n$ that commute with the face and degeneracy maps.

Now let $\Delta : Ord \to Top$ be the covariant functor such that $\Delta(\mathsf{n}):=\Delta_n$ is
the $n$-simplex and $\Delta(\mathsf{n} \stackrel{\theta}{\to} \mathsf{m}): \Delta_n \to \Delta_m$
the linear map induced by $\theta$, then we can define

\begin{definition}
The {\it geometrical realization} $|X_*|$ of the simplicial space $X$ is
$$\bigsqcup_n \Delta_n \times X_n / \sim$$
where $\sim$ is the equivalence relation $(t, X(\theta)x) \sim (\Delta(\theta)t,x)$ 
generated by all morphisms $\theta$ in $Ord$.  
\end{definition}

In the case of a groupoid $\Gg$ its nerve $\Gg_*$ form a simplicial space (see \ref{nervegroupoid})
and its realization
is what is known as the classifying space of the groupoid $B\Gg$  

\subsection{Cyclic spaces}

\begin{definition} \label{definitioncyclic}
A {\it cyclic space} $(X_*,t_*)$ consists of a simplicial space $X_*$ together with 
continuous maps $t_n: X_n \to X_n$ which satisfy the following relations with respect 
to the faces and degeneracies of $X_*$:
\begin{itemize}
\item $d_n^i \circ t_n = t_{n-1} \circ d_n^{i-1}$  and $s_n^i \circ t_n = t_{n+1} \circ s_n^{i-1}$
for $1 \leq i \leq n$ and
\item $t_n^{n+1} = {\rm id}_{X_n}$, with  $t_n^r$ the $r$-fold composition of $t_n$.
\end{itemize}
\end{definition}
The maps $t_n$ will be referred to as the cyclic structure of $X_*$. A morphism of
cyclic spaces $f_* : (X_*,t_*) \to (Y_*,t'_*)$ is a simplicial map $f_*: X_* \to Y_*$
such that $t'_n \circ f_n = f_n \circ t_n$.

\begin{rem} \label{cycliccategory}
A cyclic space can also be seen as a contravariant functor from $\Lambda$, the Connes
cyclic category \cite{Connes}, to $Top$. This category $\Lambda$ has for objects the sets $\mathsf{n}$ for $n \in \naturals$,
the same objects of $Ord$, with $Hom_\Lambda (\mathsf{n}, \mathsf{m})$ the homotopy classes of increasing continuous maps
$\phi$ of degree $1$ from $S^1$ to $S^1$ which send the $(m+1)$-st roots of unity to the $(n+1)$-st roots of unity. $Ord$ can
be seen as a subcategory of $\Lambda$.

The automorphism group of $\mathsf{n}$ in $\Lambda$ can naturally be identified with $\integer_{n+1}$, and
any morphism $\mathsf{n} \to \mathsf{m}$ can be uniquely written as $\alpha \tau^r_m$ (see \cite{Connes1}),
with $\alpha \in Hom_{Ord}(\mathsf{n},\mathsf{m})$
and $\tau_m$ the generator of $\integer_{m+1}$.

A cyclic space is a functor $(X_*,t_*) : \Lambda^{op} \to Top$; its underlying simplicial space $X_*$ is given by
restriction to $Ord^{op}$.
\end{rem}

One peculiarity of cyclic spaces is that we can endow its realization with a circle action. This
is shown in \cite{BurgheleaFiedorowicz, FiedorowiczLoday} but we will sketch its construction.

\begin{definition}
Let $X_*$ be a simplicial space. The free cyclic space on $X_*$ is the cyclic space $(S^1_* \stackrel{\sim}{\times}
X_*, t_*)$ defined by:
\begin{eqnarray*}
& & (S^1_* \stackrel{\sim}{\times} X_*)_n = \integer_{n+1} \times X_n \\
d^i_n(\tau_n^r,x) & = & \left\{ 
\begin{array}{cc}
(\tau^r_{n-1},d^{i-r}_n x) & \ if \ \ r \leq i \leq n\\
(\tau^{r-1}_{n-1},d^{n-r+i+1}_n x) & \ if \ \ 0 \leq i \leq r-1
\end{array} \right.\\
s^i_n(\tau_n^r,x) & = & \left\{ 
\begin{array}{cc}
(\tau^r_{n+1},s^{i-r}_n x) & \ if \ \ r \leq i \leq n\\
(\tau^{r+1}_{n+1},s^{n-r+i+1}_n x) & \ if \ \ 0 \leq i \leq r-1
\end{array} \right.\\
t_n(\tau^r_n,x) & =& (\tau^{r+1}_n, x).
\end{eqnarray*}
Where $0 \leq r \leq n$ and $\tau_n$ denotes the generator of $\integer_{n+1}$. We will write
$S^1_*$ for $S^1_* \stackrel{\sim}{\times} pt$.
\end{definition}

\begin{proposition} \cite[Prop 1.4]{BurgheleaFiedorowicz} \label{circleactioncyclic}
If $(X_*,t_*)$ is a cyclic space, the evaluation map $ev_*: S^1_* \stackrel{\sim}{\times} X_* \to X_*$
given by $ev_n(\tau_n^r,x)=t^r_n x$ induces a circle action $\mu:S^1 \times |X_*| \to |X_*|$.
\end{proposition}

\begin{proof} The proof can also be found in the Appendix to section 3 of \cite{Burghelea1}. We recall  it
because is used in the proof of proposition \ref{conf=inertia}.

The action $\mu$ is the composition $|ev| \circ h^{-1} \circ j^{-1}$
$$S^1 \times |X_*| = |S^1_* | \times |X_*| \stackrel{j}{\leftarrow} |S^1_* \times X_*| \stackrel{h}{\leftarrow}
|S^1_* \stackrel{\sim}{\times} X_*| \stackrel{|ev|}{\to} |X_*|$$
where the homeomorphism $h$ is given by 
\begin{eqnarray} \label{def-h}
h((\tau^r_n,x);(u_0,\dots,u_n)) = ((\tau^{n-r+1}_n,x);\tau^r_n(u_0,\dots,u_n)),
\end{eqnarray}
with $(u_0, \dots, u_n) \in \Delta_n$, $\tau_n( u_0, \dots, u_n) = (u_1, \dots, u_n, u_0)$ and
$j$ is the obvious homeomorphism from $|S^1_* \times X_*|$ to $|S^1_* | \times |X_*|$.
\end{proof}

We want to associate a topological space to a cyclic space. This is done using
homotopy colimits.

\begin{definition}
Let $\TT$ be a  small category and $C:\TT \to Top$ a functor. The { homotopy colimit}
$hocolim$ $C$ is the classifying space of the topological category $\CC$ whose space of objects is 
$\sqcup_{i \in Ob(\TT)} C(i)$, with morphisms from $x \in C(i)$ to $x' \in C(i')$ consisting
 of all morphisms $\alpha: i \to i'$ in $\TT$ such that $C(\alpha)(x) = x'$. The morphisms in $\CC$ are
topologized by identification with $\sqcup_{i,i' \in Ob(\TT)} C(i)\times Hom(i,i')$. If $X_*$ is a simplicial
space we define $\Vert X_* \Vert = hocolim \ {}_{Ord^{op}}X_*$; if $(X_*,t_*)$ is a cyclic space we define
$\Vert X_*,t_* \Vert = hocolim \ {}_{\Lambda^{op}}X_*$.
\end{definition}
For a constant functor $*: \TT \to Top$ to a point $hocolim \ * = B\TT$;  so $\Vert *,t_* \Vert = B\Lambda$.
If $(X_*,t_*)$ is a cyclic space the inclusion $Ord^{op} \subset \Lambda^{op}$ induces an inclusion
$\Vert X_* \Vert \to \Vert X_*, t_* \Vert$ and the constant map $(X_*,t_*) \to (*,t_*)$ induces a map
$\Vert X_*,t_* \Vert \to B \Lambda$. Then
\begin{theorem} \cite[Thm. 1.8]{BurgheleaFiedorowicz} \label{fibration}
Let $(X_*,t_*)$ be a cyclic space, then the sequence
$$\Vert X_* \Vert \to \Vert X_*,t_* \Vert \to B \Lambda$$
is a fibration up to homotopy, and is equivalent as a fibration to
$$ | X_* | \to ES^1\times_{S^1} |X_*| \to BS^1.$$
 \end{theorem}

In here is used the result of Connes \cite{Connes1}: $B\Lambda \cong BS^1$, and the one 
of Bousfield and Kan \cite{BousfieldKan}:
for $X_*$ a simplicial space, $\Vert X_* \Vert$ and $| X_*|$ are naturally equivalent.

\bibliographystyle{amsplain}
\bibliography{configuration}

\end{document}